\documentclass[5p]{elsarticle}

\usepackage{amsmath,amssymb}
\usepackage{mathtools}
\usepackage{tikz}
\usetikzlibrary{shapes}
\usetikzlibrary{arrows.meta}
\usetikzlibrary{arrows}
\usetikzlibrary{positioning}
\tikzset{cross/.style={cross out, draw, 
		minimum size=2*(#1-\pgflinewidth), 
		inner sep=0pt, outer sep=0pt}}

\newtheorem{lemma}{Lemma}
\newtheorem{proposition}{Proposition}
\newdefinition{remark}{Remark}
\newtheorem{assumption}{Assumption}
\newtheorem{definition}{Definition}
\newtheorem{theorem}{Theorem}
\newtheorem{example}{Example}
\newproof{proof}{Proof}

\newcommand{\norm}[1]{\left\lVert#1\right\rVert}
\newcommand\numberthis{\addtocounter{equation}{1}\tag{\theequation}}
\renewcommand{\refeq}[1]{\overset{#1}{=}}
\newcommand{\refleq}[1]{\overset{#1}{\leq}}

\newcommand{\classCost}{\mathcal{V}}
\newcommand{\algorithm}{\mathcal{A}}

\date{}

\begin{document}
	
\begin{frontmatter}
	\title{On the relation between dynamic regret and closed-loop stability\tnoteref{t1}}
	
	\tnotetext[t1]{This work was supported by the Deutsche Forschungsgemeinschaft (DFG, German Research Foundation) - 505182457.}
	
	\tnotetext[]{Published in Systems \& Control Letters, \newline DOI: https://doi.org/10.1016/j.sysconle.2023.105532}
	
	\author[1]{Marko Nonhoff\corref{cor1}}\ead{nonhoff@irt.uni-hannover.de}
	\author[1]{Matthias A. Müller}\ead{mueller@irt.uni-hannover.de}
	\cortext[cor1]{Corresponding author}
	
	\affiliation[1]{
		organization = {Leibniz University Hannover, Institute of Automatic Control},
		postcode = {30167 Hannover},
		country = {Germany}
	}
	

	\begin{abstract}
		In this work, we study the relations between bounded dynamic regret and the classical notion of asymptotic stability for the case of a priori unknown and time-varying cost functions. In particular, we show that bounded dynamic regret implies asymptotic stability of the optimal steady state for a constant cost function. For the case of an asymptotically stable closed loop, we first derive a necessary condition for achieving bounded dynamic regret. Then, given some additional assumptions on the system and the cost functions, we also provide a sufficient condition ensuring bounded dynamic regret. Our results are illustrated by examples.
	\end{abstract}
	
\end{frontmatter}

%%%%%%%%%%%%%%%%%%%%%%%%%%%%%%%%%%%%%%%%%%%%%%%%%%%%%%%%%%%%%%%%%%%%%%%%%%%%%%%%
\section{INTRODUCTION}
In recent years, there has been considerable research interest in applications of machine learning and online optimization techniques to (optimal) control problems. Along with methods and algorithms, analysis techniques that originated in the field of online learning and optimization have been applied in order to study the behavior of the closed loop consisting of a learning algorithm and a controlled system. In particular, dynamic regret
\[
	\mathcal R := \sum_{t=0}^T L_t(u_t,x_t) - L_t(\hat u_t,\hat x_t)
\]
characterizes the accumulated performance gap between the closed loop $(u_t,x_t)$ and some benchmark $(\hat u_t,\hat x_t)$ with respect to a (possibly time-varying) cost function $L_t$ to be minimized over a horizon $T$. Typically, either some part of the system's environment, most commonly either the cost functions $L_t$ or process noise affecting the system, or the system dynamics itself are assumed to be time-varying and/or a priori unknown and the benchmark $(\hat u_t,\hat x_t)$ is defined in hindsight, i.e., with full knowledge of the system dynamics and the time-varying environment. Thus, the dynamic regret~$\mathcal R$ measures the performance lost due to not knowing the environment in which the system operates a priori. It is desirable to derive an upper bound on the dynamic regret that is \textit{sublinear in $T$}, because such a bound implies that the closed loop achieves asymptotically on average the performance of the benchmark, i.e., $\lim_{T\rightarrow \infty} \mathcal R/T = 0$. In contrast to classical Lyapunov stability, which only characterizes the asymptotic behavior of the system, dynamic regret takes the closed loop's transient performance into account. Recently, regret analysis has been performed for various control algorithms, e.g., controllers minimizing dynamic regret \cite{Goel2020,Goel2021,Sabag21,Martin2022,Didier2022}, control techniques based on online convex optimization (OCO) \cite{Agarwal2019,Li2019,Li2021,Nonhoff2020,Nonhoff2021,Nonhoff2022}, and model predictive control (MPC) \cite{Wabersich2020,Dogan2021,Zhang2021,Lale2021} as well as moving horizon estimation (MHE) \cite{Gharbi2021}. However, stability results only exist for some of the above mentioned algorithms. In another closely related line of research, feedback optimization as an emerging control paradigm considers a similar setting but typically only asymptotic stability is shown instead of bounds on the dynamic regret of the closed loop \cite{Lawrence2018,Menta2018,Colombino2020,Zheng2020,Bianchin2021}.

\begin{figure}
	\begin{center}
		\begin{tikzpicture}
			\clip(-1.25,-1.5) rectangle(7.25,2.5);
			\node[align = center, text width = 3cm] at (0,2) {\underline{time-invariant} \underline{cost functions}};
			\node[align = center, text width= 3cm] at (6,2) {\underline{time-varying} \underline{cost functions}};
			\node[draw,ellipse,align=center,text width = 1.5cm] (stable) at (0,0) {asymptotic stability};
			\node[draw,ellipse,align=center, text width = 1.5cm] (regret) at (6,0) {bounded dynamic regret};
			\draw[{Stealth[length=3mm, width=2mm]}-] (stable.350) to[out=310,in=230] node[fill=white,below,yshift=-.5ex,xshift=0ex] {Theorem~\ref{thm:regret_asymp_stab}} (regret.190);
			\draw[-{Stealth[length=3mm, width=2mm]}] (stable.10) to[out=50,in=135] node[fill=white,below,yshift=-1ex,xshift=0ex,align=center,text width = 10ex] {\textit{summable} $\mathcal{KL}$ bound} node[fill=white,above,yshift = 0.5ex,text width = 22ex,align = center] {Theorem~\ref{thm:stab_regret}, Proposition~\ref{prop:necessary_condition}} (regret.170);
		\end{tikzpicture}
		\caption{Overview of our results. Summable $\mathcal{KL}$ functions are used to derive a necessary (Proposition~\ref{prop:necessary_condition}) and a sufficient (Theorem~\ref{thm:stab_regret}) condition for achieving bounded dynamic regret.}
		\label{fig:overview}
	\end{center}
		\vspace{-5ex}
\end{figure}
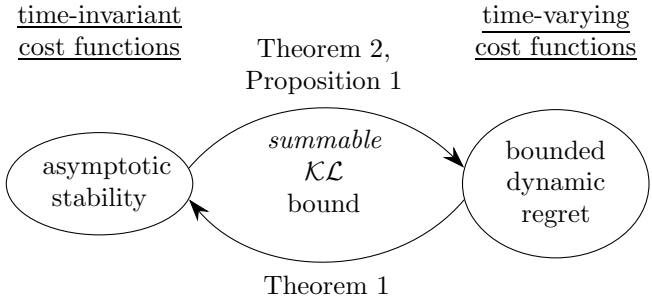

In the works mentioned above, either time-varying cost functions (due to, e.g., time-varying and a priori unknown energy prices \cite{Li2021,Nonhoff2022}), process noise (e.g., renewable energy and a priori unknown consumption in power networks \cite{Colombino2020}) or unknown system dynamics are considered. In this work, we focus on the first case, i.e., time-varying and a priori unknown cost functions, and leave process noise and learning unknown system dynamics as interesting directions for future research. As discussed above, in this setting, both dynamic regret and stability are frequently applied to study the closed loop behavior of the proposed algorithms. Therefore, it is interesting to study connections between dynamic regret and the more classical notion of stability in control theory to improve comparability and enable transferring results from different approaches that analyze one of these properties. To the authors' best knowledge, this is the first paper studying the relation between dynamic regret and stability. We show that bounded dynamic regret for time-varying cost functions implies asymptotic stability of the optimal steady state for a constant cost function under mild assumptions and that the converse implication does in general not hold. Furthermore, we derive a necessary and a sufficient condition for an asymptotically stabilizing algorithm (in the case of a constant cost function) achieving bounded dynamic regret for time-varying cost functions. Our results are illustrated in Figure~\ref{fig:overview}.

This paper is organized as follows. In Section~\ref{sec:setting}, we formally define bounded dynamic regret and asymptotic stability, and discuss the problem of interest. Section~\ref{sec:results} contains the main results of our work: We begin by proving that bounded dynamic regret implies asymptotic stability and proceed to show that the converse implication does not hold. More specifically, we provide a necessary condition for asymptotic stability implying bounded dynamic regret. Finally, we show that a similar condition is sufficient for ensuring bounded dynamic regret given some additional assumptions. Section~\ref{sec:conclusion} concludes the paper.

\textit{Notation:} The set of integer numbers and real numbers are denoted by $\mathbb{Z}$ and $\mathbb{R}$, respectively. The set of integer numbers greater than or equal to $s\in\mathbb{R}$ is $\mathbb Z_{\geq s}$. For a vector $x\in\mathbb R^n$, $\norm{x}$ denotes the Euclidean norm and for a matrix $A\in\mathbb R^{n\times m}$, $\norm{A}$ is the corresponding induced matrix norm. For a closed set $\mathcal S \subset \mathbb R^n$, $\norm{x}_\mathcal{S}$ is the distance of a point $x$ to the set $\mathcal{S}$ defined by $\norm{x}_\mathcal{S} := \inf_{y\in\mathcal{S}} \norm{x-y}$. A function $\alpha: \mathbb R_{\geq 0} \rightarrow \mathbb R_{\geq 0}$ is said to be of class $\mathcal K$ ($\alpha \in \mathcal K$) if it is continuous, strictly increasing, and $\alpha(0)=0$. If additionally $\lim_{s\rightarrow\infty}\alpha(s) = \infty$, then $\alpha\in\mathcal K_\infty$. A function $\sigma:\mathbb Z_{\geq0} \rightarrow \mathbb R_{\geq0}$ is said to be of class $\mathcal L$ if it is nonincreasing, and $\lim_{s\rightarrow\infty} \sigma(s)=0$. A function $\beta: \mathbb R_{\geq0}\times\mathbb{Z}_{\geq0}\rightarrow\mathbb R_{\geq0}$ is of class $\mathcal{KL}$ if $\beta(\cdot,t)\in\mathcal K$ for any fixed $t\in\mathbb Z_{\geq0}$ and $\beta(s,\cdot)\in\mathcal L$ for any fixed $s\in\mathbb R_{\geq0}$. A function $\beta\in\mathcal{KL}$ is called a \textit{summable} $\mathcal{KL}$-function ($\beta\in\mathcal{KL^+}$) if for any $s\in\mathbb{R}_{\geq0}$ there exists $B(s)<\infty$ such that $\sum_{t=0}^\infty \beta(s,t) \leq B(s)$.

\section{SETTING} \label{sec:setting}

In this work, we consider general discrete-time nonlinear systems of the form
\begin{equation}
	x_{t+1} = f(x_t,u_t), \label{eq:system}
\end{equation}
where $x_t\in\mathbb R^n$ and $u_t\in\mathbb R^m$ are the system states and inputs at time $t\in\mathbb Z_{\geq0}$, respectively, and \mbox{$f:\mathbb R^n \times \mathbb R^m \rightarrow \mathbb R^n$} are the system dynamics. The goal is to design an algorithm~$\algorithm$ that computes control inputs $u_t$ \textit{online} and achieves satisfactory performance with respect to the optimal control problem
\begin{align}
	\begin{split} \label{eq:OCP}
	\min_u \quad &\sum_{t=0}^T L_t(u_t,x_t) \\
	\text{s.t.} \quad &\eqref{eq:system},
	\end{split}
\end{align}
for arbitrary $T\in\mathbb Z_{\geq0}$, where the cost functions $L_t(u,x):\mathbb R^m \times \mathbb R^n \rightarrow \mathbb R$ are time-varying and a priori unknown, i.e., at each time $t$, the algorithm only has access to the previous cost functions $L_0,\dots,L_{t-1}$. Due to the a priori unknown nature of the cost functions $L_t$, we assume that the cost functions belong to a certain suitably restricted class $\classCost$, i.e., $L_t\in\classCost$ for all $t\in\mathbb Z_{\geq0}$, in order to avoid the cumulative cost $\sum_{t=0}^T L_t(u_t,x_t)$ becoming arbitrarily large. We denote the solution to \eqref{eq:OCP}, i.e., the optimal input sequence in hindsight, by $u^*(x_0) = \{u^*_t(x_0)\}_{t=0}^T$ and the corresponding state trajectory by $x^*(x_0)=\{x^*_t(x_0)\}_{t=0}^T$. The algorithm~$\algorithm$ is given by the general mapping
\[
	u_t = \algorithm(\mathcal I_t),
\]
where $\mathcal I_t = \{x_0,\dots,x_t,u_0,\dots,u_{t-1},L_0,\dots,L_{t-1},i_0\}$ describes all the available information at time $t$ and $i_0$ includes prior information, e.g., the initialization of the algorithm. In the following, we omit the dependence on $\mathcal I_t$ when it is clear from context. We make the following assumptions on the system dynamics and cost functions.
\begin{assumption} \label{assump:posdef_cost} All cost functions $L_t\in\classCost$ are positive definite with respect to a (time-varying) steady-state $(\eta_t,\theta_t)$ of~\eqref{eq:system}, i.e., there exists $(\eta_t,\theta_t)\in\mathbb R^m\times \mathbb R^n$ such that $\theta_t = f(\theta_t,\eta_t)$, $L_t(\eta_t,\theta_t) = 0$, and $L_t(u,x) > 0$ for all $(u,x)\neq (\eta_t,\theta_t)$. Moreover, there exists $\lambda \in\mathcal{K}$ such that $\lambda(\norm{(u,x)-(\eta_t,\theta_t)}) \leq L_t(u,x)$ holds for all $L_t\in\classCost$ and $(u,x)\in\mathbb{R}^m\times\mathbb R^n$.
\end{assumption}

In the following, we drop the subscript $t$ from the cost functions $L_t$ and the corresponding optimal steady state $(\eta_t,\theta_t)$ when we consider a time-invariant cost function. Assumption~\ref{assump:posdef_cost} is less restrictive than similar assumptions in the related literature. In particular, positive definiteness with respect to an optimal steady state is either assumed \cite{Nonhoff2020,Nonhoff2021} or, for cost functions that do not satisfy Assumption~\ref{assump:posdef_cost}, a strategy is chosen that tracks the optimal steady states anyway \cite{Nonhoff2022,Lawrence2018,Menta2018,Colombino2020,Zheng2020,Bianchin2021}. In the latter case, the definition of dynamic regret can be modified appropriately such that Assumption~\ref{assump:posdef_cost} holds. Moreover, the class $\mathcal K$ lower bound in Assumption~\ref{assump:posdef_cost} is commonly replaced by similar assumptions restricting the class of considered cost functions, e.g., strong convexity  \cite{Li2019,Nonhoff2022,Colombino2020}.

\begin{definition} \label{def:asymp_stability} (\cite[Definition B.8]{Rawlings2017}) The (closed, positively invariant) set $\mathcal S$ is globally asymptotically stable for a dynamical system $x_{t+1}=f(x_t,\algorithm$) if there exists a function $\beta(\cdot) \in \mathcal{KL}$ such that, for any initial state $x_0\in\mathbb R^n$,
\[
	\norm{x_t}_{\mathcal S} \leq \beta(\norm{x_0}_\mathcal{S},t) \numberthis \label{eq:def_asymp_stability}
\]
holds for all $t\in\mathbb{Z}_{\geq0}$.
\end{definition}

In addition, dynamic regret $\mathcal R$ is defined as
\begin{equation}
	\mathcal R_T(x_0,u) = \sum_{t=0}^T L_t(u_t,x_t) - L_t(\eta_t,\theta_t), \label{eq:def_regret}
\end{equation} where we omit the arguments of $\mathcal R_T(x_0,u)$ when they are clear from context. We note that in the literature, different comparator sequences than $(\eta_t,\theta_t)$ (compare Assumption~\ref{assump:posdef_cost}) are chosen regularly \cite{Agarwal2019,Li2021}. We choose $(\eta_t,\theta_t)$ because it is the hardest benchmark, i.e., dynamic regret with respect to any other sequence is lower than or equal to the dynamic regret with respect to $(\eta_t,\theta_t)$ by Assumption~\ref{assump:posdef_cost}. Furthermore, bounded dynamic regret of the Algorithm~$\algorithm$ with respect to the benchmark $(\eta_t,\theta_t)$ is equivalent to bounded dynamic regret of the Algorithm~$\algorithm$ with respect to any other benchmark, if this benchmark achieves bounded regret itself, compare Remark~\ref{rem:different_benchmark} below. Then, the goal is to achieve a \textit{sublinear} (in $T$) upper bound of the dynamic regret, because such a bound implies that the algorithm in a closed loop with the system achieves asymptotically on average a cost that is no worse than the optimal cost. However, \cite{Li2019} indicates that the best achievable bound on the dynamic regret is linear in the path length, where the path length $PL$ is given by
\[
	PL = \sum_{t=1}^T \norm{\eta_{t}-\eta_{t-1}} + \norm{\theta_{t}-\theta_{t-1}}. \numberthis \label{eq:path_length}
\]
The path length can be seen as an indicator of the variation of the cost functions $L_t$. Thus, we define bounded dynamic regret as follows.

\begin{definition} \label{def:bounded_regret} An algorithm~$\algorithm$ achieves bounded dynamic regret if, for any initial state $x_0\in\mathbb R^n$, there exist constants $C_\eta,C_\theta,C \geq 0$, independent of $T$, such that
\[
	\mathcal R_T(x_0,u) \leq C_\eta \sum_{t=1}^T \norm{\eta_t - \eta_{t-1}} + C_\theta \sum_{t=1}^T \norm{\theta_t - \theta_{t-1}} + C
\]
holds for any sequence of cost functions $L_0, \dots, L_T\in\classCost$.
\end{definition}

\begin{remark} \label{rem:different_benchmark} (Other benchmark trajectories) As discussed above, other benchmarks than the optimal steady states $(\eta_t,\theta_t)$ are commonly applied in the related literature. Let
\[
	\hat{\mathcal R}_T := \sum_{t=0}^T L_t(u_t,x_t) - L_t(\hat u_t, \hat x_t)
\]
be the dynamic regret with respect to the benchmark $(\hat u,\hat x)$. Moreover, assume that the benchmark trajectory $(\hat u,\hat x)$ achieves bounded dynamic regret according to Definition~\ref{def:bounded_regret}, i.e., there exist $C_\eta^*, C_\theta^*, C^* \geq 0$ such that
\begin{align*}
	 &\sum_{t=0}^T L_t(\hat u_t, \hat x_t) - L_t(\eta_t, \theta_t) \\
	 \leq &C^*_\eta \sum_{t=1}^T \norm{\eta_t - \eta_{t-1}} + C^*_\theta \sum_{t=1}^T \norm{\theta_t - \theta_{t-1}} + C^*
\end{align*}
holds for all $T\in\mathbb Z_{\geq1}$ and $L_0,\dots,L_T\in\classCost$.	For example, when considering the optimal trajectory in hindsight as a benchmark, i.e., $(\hat u, \hat x) = (u^*(x_0),x^*(x_0))$, similar assumptions on the optimal achievable performance are common in the related literature, compare, e.g., \cite{Agarwal2019,Li2021}. Then, the algorithm~$\algorithm$ achieves bounded dynamic regret (i.e., $\mathcal R_T \leq C_\eta \sum_{t=1}^T \norm{\eta_t - \eta_{t-1}} + C_\theta \sum_{t=1}^T \norm{\theta_t - \theta_{t-1}} + C$) if and only if it achieves bounded regret with respect to the benchmark trajectory $(\hat u,\hat x)$ (i.e., $\hat{\mathcal R}_T \leq \hat C_\eta \sum_{t=1}^T \norm{\eta_t - \eta_{t-1}} + \hat C_\theta \sum_{t=1}^T \norm{\theta_t - \theta_{t-1}} + \hat C$). This claim can be proven by noting that i) $\hat{\mathcal R}_T \leq \mathcal R_T$ due to Assumption~\ref{assump:posdef_cost} and ii) $\mathcal R_T = \hat{\mathcal R}_T + \sum_{t=0}^T L_t(\hat u_t,\hat x_t) - L_t(\eta_t,\theta_t)$, which implies the result using the above assumption that the benchmark's dynamic regret is bounded.
\end{remark}

\section{MAIN RESULTS} \label{sec:results}
In this section, we derive conditions under which an algorithm~$\algorithm$ that achieves bounded dynamic regret according to Definition~\ref{def:bounded_regret} (for time-varying cost functions) also asymptotically stabilizes the optimal steady state $(\eta,\theta)$ (for a time-invariant cost function) and vice versa. Note that, in the case that the cost functions are constant, i.e., $L_t(u,x) = L(u,x)$ for $t\in\mathbb Z_{\geq0}$, convergence to the optimal steady state, i.e., $\lim_{t\rightarrow\infty} (u_t,x_t) = (\eta,\theta)$, immediately follows from the definition of bounded regret in Definition~\ref{def:bounded_regret} and Assumption~\ref{assump:posdef_cost}.

\subsection{Bounded regret implies stability}
First, we analyze algorithms that are already known to achieve bounded dynamic regret. In order to also establish asymptotic stability, consider the following assumptions.

\begin{assumption}\label{assump:algorithm_state_space_rep} The algorithm $\algorithm$ admits a state-space representation, i.e., it can be written as
\begin{align} \label{eq:algorithm_state_space}
	\begin{split}
	x_{t+1}^{\algorithm} &= f^{\algorithm}_{L_{t-1}}(x_t^{\algorithm},x_t) \\
	u_t &= h^\algorithm_{L_{t-1}}(x^{\algorithm}_{t}, x_t),
	\end{split}
\end{align}
where $x^{\algorithm} \in \mathbb R^{n_\algorithm}$, $f^{\algorithm}_{L_{t-1}}: \mathbb R^{n_{\algorithm}} \times \mathbb R^n \rightarrow \mathbb R^{n_\algorithm}$, and $h^\algorithm_{L_{t-1}}: \mathbb R^{n_\algorithm} \times \mathbb R^n \rightarrow \mathbb R^m$.
\end{assumption}

Assumption~\ref{assump:algorithm_state_space_rep} restricts the class of considered algorithms. As standard in OCO, we assume that the algorithm only depends on the previous cost function $L_{t-1}$ because the current cost function $L_t$ is unknown \cite{Li2019,Nonhoff2022}. However, this assumption is not too restrictive and includes a wide range of nonlinear algorithms that admit a state-space realization. Moreover, we are still able to give guarantees for algorithms that do not satisfy Assumption~\ref{assump:algorithm_state_space_rep} as discussed in Remark~\ref{rem:asymp_stab_general_algo} below.

Furthermore, we denote by $\mathcal S_L \subseteq \mathbb R^{n_\mathcal{A}}$ the largest positively invariant set of states of the algorithm~$\algorithm$ that produce the control input $\eta$ for \mbox{$x = \theta$} and a (constant) cost function~$L\in\classCost$, i.e., 
	\begin{align*}
		\mathcal S_L := \{ x^{\algorithm} ~ &| ~ \chi^\algorithm_0 = x^\algorithm, \\
		&\eta = h_L^\algorithm(\chi^\algorithm_\tau,\theta),~ \chi^\algorithm_{\tau+1} = f^{\algorithm}_L(\chi^{\algorithm}_\tau,\theta) ~\forall\tau\in\mathbb Z_{\geq0} \},
	\end{align*}
and we define 
\[
	\mathcal S^\theta_L := \{ (x,x^\algorithm) ~|~ x=\theta, x^\algorithm \in S_L \}.
\]
Then, we require the constant $C$ in Definition~\ref{def:bounded_regret} to depend on the initial states of the closed loop as follows.

\begin{assumption} \label{assump:class_k_constant}
	The Algorithm~$\algorithm$ and the set of cost functions~$\classCost$ are such that the constant $C$ in Definition~\ref{def:bounded_regret} satisfies $C = C(x_0,x_0^{\mathcal A})$ and for any constant cost function $L\in\classCost$, there exists $\bar\alpha \in \mathcal K_\infty$ such that
	\[
		C(x_0,x_0^{\algorithm}) \leq \bar\alpha(\norm{(x_0,x_0^{\algorithm})}_{\mathcal S^\theta_{L}})
	\]
	holds for any $(x_0,x_0^\algorithm)\in\mathbb R^n\times\mathbb R^{n_\algorithm}$.	
\end{assumption}

Assumption~\ref{assump:class_k_constant} requires the constant $C$ in Definition~\ref{def:bounded_regret} to solely depend on the initial state $x_0$ and the algorithm's initialization. It does, however, not require the closed loop to stay at $(\eta,\theta)$ for $(u_0,x_0) = (\eta,\theta)$ if $x_0^\algorithm \notin \mathcal S_{L}$. Finally, we assume a lower bound on the dynamic regret. Therein, with a slight abuse of notation, we write $\mathcal R_{T}(x_0,x_0^\algorithm)$ for $\mathcal R_{T}(x,u)$, because the input sequence $u$ is uniquely defined by the initial states $(x_0,x_0^\algorithm)$ by Assumption~\ref{assump:algorithm_state_space_rep}.

\begin{assumption} \label{assump:regret_lower_bound} 
		The algorithm~$\algorithm$ and the set of cost functions~$\classCost$ are such that for any constant cost function $L\in\classCost$, there exists a function $\underline \alpha \in \mathcal K_\infty$ that satisfies
	\[
		\underline\alpha(\norm{(x_0,x_0^\algorithm)}_{\mathcal S_{L}^\theta}) \leq \lim_{T\rightarrow\infty} \mathcal R_T(x_0,x_0^\algorithm).
	\]
	for any $(x_0,x_0^\algorithm)\in\mathbb R^n\times\mathbb R^{n_\algorithm}$.
\end{assumption}

Since we do not pose any continuity assumptions on the algorithm dynamics $f^\algorithm_L$ and $h^\algorithm_L$, we require Assumption~\ref{assump:regret_lower_bound} in order to avoid discontinuities such that, for a constant cost function $L\in\classCost$, $\lim_{t\rightarrow\infty} (u_t,x_t) = (\eta,\theta)$, but $\lim_{t\rightarrow\infty} (x_t,x_t^\algorithm) \notin S_L^\theta$. For example, Assumption~\ref{assump:regret_lower_bound} is satisfied on a compact set (i.e, for all $(x_0,x_0^\algorithm)\in\mathcal X_0^\algorithm$, where $\mathcal X_0^\algorithm\subseteq \mathbb R^{n+n_\algorithm}$ is a compact set) if the algorithm dynamics $f^\algorithm_L$ and $h^\algorithm_L$, the system dynamics $f$, and the cost function $L$ are continuous. Then, using continuity of the dynamics and the cost function together with Assumption~\ref{assump:class_k_constant}, one can show by applying the uniform limit theorem \cite[Theorem 21.6]{Munkres2000} that $\lim_{T\rightarrow\infty} \mathcal R_T(x_0,x_0^\algorithm)$ is a continuous function, which yields the desired bound by Assumption~\ref{assump:posdef_cost}.

Note that Assumptions~\ref{assump:algorithm_state_space_rep}-\ref{assump:regret_lower_bound} are satisfied for, e.g, the algorithms proposed in \cite{Nonhoff2020,Nonhoff2022} by choosing the predicted input sequence and the estimated optimal input and steady state, respectively, as the controller states $x^\algorithm$. Our first main result states that any algorithm~$\algorithm$ that achieves bounded regret as specified by Definition~\ref{def:bounded_regret} asymptotically stabilizes the closed loop.

\begin{theorem} \label{thm:regret_asymp_stab} 
	Suppose Assumptions \ref{assump:posdef_cost}-\ref{assump:regret_lower_bound} are satisfied. For any constant cost function $L\in\classCost$, the set $\mathcal S_L^\theta$ is globally asymptotically stable for the extended state $(x,x^\algorithm)$ with respect to the closed loop dynamics given by~\eqref{eq:system} and~\eqref{eq:algorithm_state_space} if the algorithm achieves bounded dynamic regret for time-varying cost functions.
\end{theorem}

\begin{proof}
	Fix any cost function $L\in\classCost$ and assume that the cost function is time invariant, i.e., $L_t(u,x) = L(u,x)$ and $(\eta_t,\theta_t) = (\eta,\theta)$ for all $t\in\mathbb{Z}_{\geq0}$. We begin by choosing
	\[
		V(x,x^\algorithm) = \lim_{T\rightarrow \infty} \mathcal R_{T}(x,x^\algorithm), \numberthis \label{eq:def_V}
	\]
	as a Lyapunov function candidate. Note that the above limit exists since $(\eta_t,\theta_t) = (\eta,\theta)$ holds for all $t\in\mathbb{Z}_{\geq0}$ and the algorithm achieves bounded dynamic regret (compare Definition~\ref{def:bounded_regret}), i.e.,
	\[
		\lim_{T \rightarrow \infty} \mathcal R_{T}(x,x^\algorithm) \leq C. \numberthis \label{eq:bounded_regret_limit}
	\]
	Moreover, by Assumptions~\ref{assump:class_k_constant} and~\ref{assump:regret_lower_bound} there exist functions $\underline \alpha, \bar \alpha \in \mathcal K_\infty$ such that 
	\begin{equation}
		\underline \alpha(\norm{(x,x^\algorithm)}_{\mathcal S^\theta_L}) \leq V(x,x^\algorithm) \refleq{\eqref{eq:bounded_regret_limit}} \bar\alpha(\norm{(x,x^{\algorithm})}_{S^\theta_L}) . \label{eq:bounded_lyap_function}
	\end{equation}
	In the following, we define $V_t = V(x_t,x^\algorithm_t)$. Then, we get
	\begin{align*}
		&V_{t+1} - V_{t} = \lim_{T\rightarrow \infty} \mathcal R_{T}(x_{t+1},x^{\algorithm}_{t+1}) - \lim_{T\rightarrow \infty} \mathcal R_{T}(x_t,x^{\algorithm}_t) \\
		\refeq{\eqref{eq:bounded_regret_limit}} &\lim_{T\rightarrow \infty} \left( \mathcal R_{T}(x_{t+1},x^{\algorithm}_{t+1}) - \mathcal R_{T}(x_t,x^{\algorithm}_t) \right) \\
		\begin{split}
		\refeq{\eqref{eq:def_regret}} &\lim_{T\rightarrow \infty} \left( \sum_{\tau = t+1}^{T+t+1} \left( L(h_L^\algorithm(x^{\algorithm}_{\tau},x_{\tau}),x_{\tau}) - L(\eta,\theta) \right) \right. \\ & \left.- \sum_{\tau=t}^{T+t} \left( L(h_L^\algorithm(x^{\algorithm}_{\tau},x_{\tau}),x_{\tau}) - L(\eta,\theta) \right) \right)
		\end{split} \\
		\begin{split}
		= &\lim_{T\rightarrow \infty} \Big( L(h_L^\algorithm(x^\algorithm_{T+t+1},x_{T+t+1}),x_{T+t+1}) \\ &\qquad - L(h^\algorithm(x^\algorithm_{t},x_{t}),x_{t}) \Big).
		\end{split}
	\end{align*}
	As discussed above, bounded dynamic regret implies convergence, i.e., $\lim_{T\rightarrow\infty} L(h_L^\algorithm(x^\algorithm_{T},x_{T}),x_{T}) = L(\eta,\theta)$. Thus, we get by Assumption~\ref{assump:posdef_cost}
	\begin{align*}
		V_{t+1} - V_{t} &= - \Big( L(h_L^\algorithm(x^\algorithm_{t},x_{t}),x_{t}) -L(\eta,\theta) \Big) \leq 0.
	\end{align*}
	Together with~\eqref{eq:bounded_lyap_function}, this implies stability of the set $\mathcal S_L^\theta$ by standard Lyapunov arguments. It remains to show convergence of the extended state to $\mathcal S_L^\theta$. We have $\lim_{t\rightarrow\infty} V_t = 0$, because the dynamic regret $V_0 = \lim_{T\rightarrow\infty} \mathcal R_T(x_0,x_0^\algorithm)$ is bounded by a finite constant $C$ by Definition~\ref{def:bounded_regret}, which implies $\lim_{t\rightarrow\infty} \norm{(x_t,x_t^\algorithm)}_{\mathcal S_L^\theta} = 0$ by Assumption~\ref{assump:regret_lower_bound}. \hfill $\square$
	\end{proof}
\begin{remark} \label{rem:asymp_stab_general_algo} (Algorithms that admit a state-space representation) If Assumption~\ref{assump:algorithm_state_space_rep} is not satisfied, it is still possible to show asymptotic stability of the states of the system: To this end, instead of Assumption~\ref{assump:class_k_constant}, we assume \linebreak $C = C(x_0) \leq \bar\alpha(\norm{x_0-\theta_0})$ for some $\bar\alpha \in \mathcal K_\infty$. Note that this is a considerably stronger assumption, because Assumption~\ref{assump:class_k_constant} allows the states $x_t$ to initially diverge from $\theta_0$, even if the initial state $x_0$ is close to $ \theta_0$. Thus, if the algorithm does not admit a state-space representation, we have to assume that it is initialized correctly, which is a restrictive assumption since the pair $(\eta_0,\theta_0)$ is a priori unknown as discussed above. Moreover, instead of Assumption~\ref{assump:regret_lower_bound}, we require $\underline \alpha(\norm{x_0-\theta_0}) \leq \lim_{T\rightarrow\infty} R_T(x_0,u)$ for some $\underline\alpha\in\mathcal K_\infty$, which is satisfied if the function $\lambda$ in Assumption~\ref{assump:posdef_cost} is of class $\mathcal K_\infty$ (i.e., the cost function $L$ is radially unbounded) because $\lim_{T\rightarrow\infty} R_T(x_0,u) \geq L(u_0,x_0) \geq \lambda(\norm{(u_0,x_0)-(\eta_0,\theta_0)}) \geq \lambda(\norm{x_0-\theta_0})$. Then, similar arguments to those in the proof of Theorem~\ref{thm:regret_asymp_stab} can be applied to show that the optimal steady state $\theta$ is asymptotically stable. Moreover, for constant cost functions $L(u,x)$, the input to the system $u_t$ still converges to $\eta$ as discussed above and the deviation of $u_t$ from $\eta$ is bounded, since we have $R_T(x_0,u) \leq C$ by~\eqref{def:bounded_regret} and $R_T(x_0,u)$ is a positive definite function with respect to $\norm{u-\eta}$ by Assumption~\ref{assump:posdef_cost}.
\end{remark}

\begin{remark} (Other benchmark trajectories)
	As discussed in Remark~\ref{rem:different_benchmark}, Theorem~\ref{thm:regret_asymp_stab} also holds for different benchmark trajectories in the definition of dynamic regret, if the benchmark achieves bounded dynamic regret according to Definition~\ref{def:bounded_regret}. If this is not the case, Theorem~\ref{thm:regret_asymp_stab} still holds true if, for any constant cost function $L\in\classCost$, i) there exists an optimal solution $u^*(x_0)$ to the OCP~\eqref{eq:OCP} for every $x_0\in\mathbb R^n$ and any $T\in\mathbb Z_{\geq0}\cup\infty$, ii) the benchmark is chosen to be the optimal trajectory in hindsight $(u^*(x_0),x^*(x_0))$, iii) Assumptions~\ref{assump:posdef_cost}-\ref{assump:regret_lower_bound} are satisfied (for the modified definition of dynamic regret), and iv) the optimal trajectory converges to the optimal steady state for every $x_0\in\mathbb R^n$, i.e., $\lim_{t\rightarrow\infty} (u^*_t(x_0),x^*_t(x_0)) = (\eta,\theta)$. To show this, we define as a Lyapunov function candidate 
	\[
		V(x) = \lim_{T \rightarrow \infty} \sum_{\tau=t}^{T+t} L(u_\tau,x_\tau) - L(u^*_\tau(x_t),x^*_\tau(x_t)).
	\]	
	The Lyapunov function $V(x)$ is then upper (lower) bounded by $\bar \alpha\,(\underline\alpha)\in\mathcal K_\infty$ due to Assumptions~\ref{assump:class_k_constant} and~\ref{assump:regret_lower_bound}, respectively. Moreover, we have $V(x_{t+1}) - V(x_t)\leq0$ because
	\begin{align*}
	\begin{split}
		\lim_{T \rightarrow \infty} &\left( \sum_{\tau=t}^{T+t} L(u^*_\tau(x_t),x^*_\tau(x_t)) \right. \\
		\quad &\left. - L(u_t,x_t) - \sum_{\tau=t+1}^{T+t} L(u^*_{\tau}(x_{t+1}),x^*_{\tau}(x_{t+1})) \right) \leq 0
	\end{split}
	\end{align*}
	by optimality of $(u^*(x_t),x^*(x_t))$. Thus, we can conclude stability of the optimal steady state $\theta$ by standard Lyapunov arguments. Convergence follows by similar arguments as in the proof of Theorem~\ref{thm:regret_asymp_stab}.
\end{remark}

\subsection{A necessary condition for bounded regret}

Theorem~\ref{thm:regret_asymp_stab} and Remark~\ref{rem:asymp_stab_general_algo} show that, under certain conditions, bounded dynamic regret implies asymptotic stability. However, the converse implication does in general not hold, i.e., asymptotic stability with respect to an invariant cost function does not imply bounded regret for time-varying cost functions. This is due to the fact that asymptotic stability only characterizes the asymptotic behavior of the closed loop, whereas bounded regret requires a certain rate of convergence, which is illustrated by the next example.

\begin{example} \label{cex:stable->regret}
	Consider the discrete-time integrator \linebreak $x_{t+1} = x_t + u_t$. We design an algorithm given by 
	\[
		u_t = -\frac{1}{\tau+1} x_t + \frac{1}{\tau+1}\theta_{t-1}, 
	\]
	where $\theta_{-1} = x_0$ and $\tau \in \mathbb Z_{\geq0}$ is set to $\tau = 0$ at time $t=0$ and increased by one at every time step. When a changed cost function is revealed, we reset $\tau$ to one. If the minimizer $\theta$ is constant, it can be verified by repeatedly inserting the system dynamics that the closed loop achieves $\norm{x_t - \theta} = \frac{1}{t} \norm{x_0 - \theta}$ for $t\in\mathbb Z_{\geq1}$. Therefore, the closed loop is globally asymptotically stable for any time-invariant cost function $L\in\classCost$ according to Definition~\ref{def:asymp_stability} with $\mathcal S = \{\theta\}$ and $\beta(s,t) = \frac{s}{t}$ if $t>1$ and $\beta(s,t) = s$ if $t=0$. However, by choosing the cost function $L(u,x) = \norm{x - \theta}+\norm{u}$ and neglecting the term $\norm{u}$, we get a lower bound for the regret
	\[
		\mathcal R = \sum_{t=0}^T L(u_t,x_t) \geq \sum_{t=0}^T \norm{x_t-\theta} = \norm{x_0-\theta} \left( 1 + \sum_{t=1}^T \frac{1}{t} \right)
	\]
	which cannot be bounded by a constant independent of $T$ because the sum $\sum_{t=1}^T \frac{1}{t}$ diverges.
\end{example}

We can generalize Example~\ref{cex:stable->regret} in order to get a necessary condition for achieving bounded regret. For this, we require the following assumption.

\begin{assumption} \label{assump:bounded_convergence}
	The algorithm~$\algorithm$ and the set of cost functions~$\classCost$ are such that, for any constant cost function $L$, there exists $\underline T \in \mathbb I_{\geq0}$ and $\mu \geq 1$ such that, for any $s\geq0$, there exists a closed-loop trajectory $\bar x =  \{ \bar x_{t} \}_{t=0}^\infty$ which satisfies
	\[
		\mu \lim_{T \rightarrow \infty} \sum_{t=\underline T}^T \norm{\bar x_{t} - \theta} \geq \lim_{T \rightarrow \infty} \sum_{t=\underline T}^T \phi(s,t),
	\]
	where 
	\begin{equation} \label{eq:def_phi}
		\phi(s,t) = \sup_{x_0,\tau} \norm{x_\tau - \theta} ~\text{s.t.}~ \norm{x_0-\theta}\leq s,~\tau \geq t.
	\end{equation}
\end{assumption}

Assumption~\ref{assump:bounded_convergence} ensures that the accumulated distance to the optimal steady state $\theta$ of any trajectory that starts in a bounded set around $\theta$ can be upper bounded by that of another trajectory potentially starting farther away from $\theta$. Thus, Assumption~\ref{assump:bounded_convergence} excludes pathological cases where the convergence rate of the closed loop is arbitrarily slow, e.g., when the initial state is close to $\theta$ (see below for an example). We use Assumption~\ref{assump:bounded_convergence} to show that, if every trajectory of a closed-loop system $x_{t+1} = f(x_t,\algorithm)$ is summable, then there exists a summable $\mathcal{KL}$ function $\beta\in\mathcal{KL}^+$ in Definition~\ref{def:asymp_stability}.

\begin{lemma} \label{lem:summable_KL_function}
	Let Assumptions~\ref{assump:posdef_cost} and~\ref{assump:bounded_convergence} hold and consider a constant cost function $L\in\classCost$. Assume that the set~$\mathcal S=\{\theta\}$ is asymptotically stable with respect to the closed-loop dynamics $x_{t+1} = f(x_t,\algorithm)$. Then, there exists a function $\beta\in\mathcal{KL^+}$ such that 
	\[
		\norm{x_t-\theta} \leq \beta(\norm{x_0-\theta},t) \numberthis \label{eq:proposition}
	\]
	holds for all $t\in\mathbb Z_{\geq0}$ and $x_0\in\mathbb R^n$ if and only if for any $x_0\in\mathbb R^n$ there exists a constant $D(\norm{x_0-\theta})$ such that
	\[
		\lim_{T\rightarrow\infty} \sum_{t=0}^T \norm{x_t-\theta} \leq D(\norm{x_0-\theta})<\infty.
	\]
\end{lemma}

The proof of Lemma~\ref{lem:summable_KL_function} is given in the Appendix.

If Assumption~\ref{assump:bounded_convergence} is not satisfied, Lemma~\ref{lem:summable_KL_function} does in general not hold. For example, a system that satisfies $\norm{x_t-\theta} = s \left( \frac{1}{1+s} \right)^t$ for a constant cost function $L$, where $s=\norm{x_0-\theta}$, is asymptotically (exponentially) stable for any $x_0\in\mathbb R^n$. Moreover, $\sum_{t=0}^\infty \norm{x_t-\theta} = \sum_{t=0}^\infty s\left(\frac{1}{1+s}\right)^t = s+1 = D(s)$ holds for any $s\geq0$, i.e., the latter condition of Lemma~\ref{lem:summable_KL_function} is satisfied. However, it can be shown that there does not exist a $\beta\in\mathcal{KL}^+$ (but only $\beta \in \mathcal{KL}$) that satisfies~\eqref{eq:proposition}. 

Next, we state a necessary condition for an asymptotically stabilizing algorithm achieving bounded regret.

\begin{proposition} \label{prop:necessary_condition}
	Let Assumptions~\ref{assump:posdef_cost} and~\ref{assump:bounded_convergence} hold. Assume that, for any constant cost function~$L\in\classCost$, the set~$\mathcal S=\{\theta\}$ is asymptotically stable with respect to the closed-loop dynamics $x_{t+1} = f(x_t,\algorithm)$. If~$\algorithm$ achieves bounded dynamic regret, then it holds that for every constant cost function~$L\in\classCost$ that satisfies $L(u,x) \geq \bar\mu\norm{x-\theta}$ for any $(u,x)\in\mathbb R^m\times\mathbb R^n$ and some $\bar\mu>0$, there exists $\beta\in\mathcal{KL^+}$ such that~\eqref{eq:proposition} is satisfied for all $t\in\mathbb Z_{\geq0}$ and $x_0 \in \mathbb R^n$.
\end{proposition}

\begin{proof}
	We prove this proposition by contraposition, i.e., if~\eqref{eq:proposition} does not hold for any $\beta \in \mathcal{KL^+}$ and a constant cost function~$L\in\classCost$ that satisfies the requirements of Proposition~\ref{prop:necessary_condition}, then $\algorithm$ cannot achieve bounded regret. Fix any constant cost function $L\in\classCost$ and $\bar\mu>0$ such that $L(u,x)\geq \bar\mu\norm{x-\theta}$ for all $(u,x)\in\mathbb R^m\times\mathbb R^n$, and assume that there does not exist any $\beta\in\mathcal{KL^+}$ that satisfies~\eqref{eq:proposition}. By Lemma~\ref{lem:summable_KL_function}, there exists an initial state $x_{0}$ such that
	\[
		\lim_{T\rightarrow \infty} \sum_{t=0}^T \norm{x_{t} - \theta} = \infty.
	\]
	The dynamic regret is then lower bounded by
	\[
		\lim_{T\rightarrow\infty} \mathcal R_T(x_0,u) \geq \lim_{T\rightarrow\infty} \sum_{t=0}^T \bar \mu \norm{x_{t} - \theta} = \infty,
	\]
	i.e., it cannot be bounded independent of $T$ since the sum diverges. Therefore, the algorithm cannot achieve bounded dynamic regret. \hfill $\square$
\end{proof}

Proposition~\ref{prop:necessary_condition} states that a necessary condition for \linebreak achieving bounded dynamic regret is existence of \textit{summable} $\mathcal{KL}$ functions $\beta$ satisfying~\eqref{eq:proposition} for a certain class of constant cost functions. We note that summable $\mathcal{KL}$ functions are, e.g., also used in the context of economic MPC \cite{Grune2013} and MHE \cite{Knuefer2021}.

\subsection{A sufficient condition for bounded regret}

Finally, we derive a sufficient condition for an asymptotically stabilizing algorithm to achieve bounded regret. To this end, we let $N$ be the number of time steps in the time interval $[0,T]$ where the cost function $L_t$ switches (including $t=0$) and let $t_i$, $i\in[0,N-1]$, be the time instances where the cost function $L_t$ changes, i.e., $t_0 = 0$ and
\begin{align*}
	L_t(u,x) \neq L_{t-1}(u,x) &\text{ if and only if } t = t_i,~i\in[1,N-1] \\
	L_t(u,x) = L_{t_i}(u,x) &\text{ if } t \in [t_i,t_{i+1}-1].
\end{align*}
To simplify notation, we additionally define $t_N = T$. We require that the cost functions are not allowed to change too frequently and restrict the function $\beta\in\mathcal{KL}$ in Definition~\ref{def:asymp_stability} to be linear in $\norm{x_0-\theta}$.

\begin{assumption} \label{assump:avg_dwell_time}
	Denote by $\mathcal N(\tau_1,\tau_2)$ the number of times the cost function switches in the interval $[\tau_1,\tau_2]\subseteq \mathbb Z_{[0,T]}$. There exist constants $N_0 \in \mathbb I_{\geq0}$ and $\varphi > 0$ such that 
	\[
	\mathcal N(\tau_1,\tau_2) < N_0 + \frac{1}{\varphi}(\tau_2-\tau_1). \numberthis \label{eq:avg_dwell_time}
	\]
\end{assumption}

\begin{assumption} \label{assump:state_sequence}
	The Algorithm~$\algorithm$ and the set of cost functions~$\classCost$ are such that there exist $\sigma\in\mathcal L$, $\sigma(0) = 1$, and a constant $k \geq 1$ such that the algorithm~$\algorithm$ in closed loop with system~\eqref{eq:system} satisfies
	\[
	\norm{x_t - \theta} \leq k\norm{x_0-\theta}\sigma(t)
	\]
	for all $t\in\mathbb Z_{\geq 0}$ and any constant cost function $L\in\classCost$.
\end{assumption}

Assumption~\ref{assump:avg_dwell_time} is commonly applied in the literature on switched systems to guarantee stability (compare, e.g., \cite{Liberzon2003}). Therein, the constants $N_0$ and $\varphi$ are referred to as \textit{chatter bound} and \textit{average dwell-time}, respectively. In our setting, such an assumption is necessary because, if the cost function $L_t$ was allowed to change at every time step, then Assumption~\ref{assump:state_sequence} would not be sufficient to guarantee convergence since the constant~$k$ in Assumption~\ref{assump:state_sequence} is larger than one. However, it would be reasonable to assume that the previous bound on $\norm{x_t-\theta_{t_i}}$ does not increase much if $\theta_{t_{i+1}}$ is close to $\theta_{t_i}$, i.e., the new upper bound on $\norm{x_t - \theta_{t_i+1}}$ from Assumption~\ref{assump:state_sequence} specifying the rate of convergence during the interval $[t_{i+1},t_{i+2}-1]$ should depend on $\norm{\theta_{t_{i+1}}-\theta_{t_{i}}}$. We conjecture that such an assumption could render Assumption~\ref{assump:avg_dwell_time} unnecessary. We leave this problem as an interesting direction for future research.

Moreover, note that we require a bound in Assumption~\ref{assump:state_sequence} which is \textit{uniform} in $L$. In order to obtain such a bound, it is necessary to restrict the class of considered cost functions. Note that this is analogous to similar assumptions that are typically taken in the literature, which also require some uniform properties for all considered cost functions, such as, e.g., (strong) convexity (compare, e.g., \cite{Agarwal2019,Li2019,Li2021,Nonhoff2022, Colombino2020}). 

Finally, assume that the initial state and the optimal steady states are contained in compact sets $\mathcal X_0, \Theta \subset \mathbb R^n$, i.e., $x_0\in\mathcal X_0$ and $\theta_t \in \Theta$ for all $\mathbb Z_{\geq0}$, and that $\varphi$ in Assumption~\ref{assump:avg_dwell_time} is large enough. Then, there exists a compact set $\mathcal X \in \mathbb R^n$ such that $x_t \in \mathcal X$ for all $t$. Moreover, for any function $\beta\in\mathcal{KL}$ there exist $\kappa\in\mathcal K_\infty$ and $\sigma\in\mathcal L$ such that $\beta(s,t)\leq \kappa(s)\sigma(t)$ \cite[Lemma~8]{Sontag1998}. Therefore, Assumption~\ref{assump:state_sequence} is satisfied if a \textit{uniform} bound (in $L$) can be found such that $\kappa(s)$ is Lipschitz at the origin.

In order to make use of Assumption~\ref{assump:avg_dwell_time}, we need to ensure that the average dwell time $\varphi$ is sufficiently large.

\begin{lemma} \label{lem:bounded_product}
	Let Assumptions \ref{assump:avg_dwell_time} and~\ref{assump:state_sequence} be satisfied. For every $N_0\in\mathbb{Z}_{\geq0}$, there exists $\underline \varphi > 0$, $P>0$, such that for every $\varphi \geq \underline\varphi$ and any $N' \in \mathbb Z_{>N_0}$, it holds that
	\[
	\sum_{i=0}^{N'}\prod_{j=1}^{i} k(\Delta \sigma)_j \leq P, \numberthis \label{eq:bound_sum_prod}
	\]
	where $(\Delta \sigma)_j = \sigma(t_j-t_{j-1})$.
\end{lemma}

The proof is given in the Appendix. 

In addition to asymptotic stability, bounded dynamic regret requires bounded inputs due to their influence on the cost functions.

\begin{assumption} \label{assump:bounded_inputs}
	The algorithm~$\algorithm$ and the set of cost functions~$\classCost$ are such that the inputs $u_t$ satisfy
	\[
	\begin{split} 
		&\norm{u_t-\eta_{t-1}} \leq k_u \norm{u_{t-1}-\eta_{t-1}} \\&\qquad+ k_x \norm{x_t-\theta_{t-1}} + k_\zeta \norm{\zeta_{t-1} - \zeta_{t-2}} 
	\end{split} \numberthis \label{eq:bounded_inputs}
	\]
	for all $t\in\mathbb Z_{\geq0}$, where $\zeta_t = \begin{bmatrix} \theta_t^\top & \eta_t^\top \end{bmatrix}^\top$, $\zeta_{t} = \zeta_0$ if $t\in\mathbb Z_{<0}$, $u_{-1} = \eta_0$, $k_u \in [0,1)$ and $k_x,k_\zeta \geq0$.
\end{assumption}

In time intervals during which the algorithm receives a constant cost function $L_{t_i}$, i.e., $t\in[t_i+2,t_{i+1}]$, Assumption~\ref{assump:bounded_inputs} requires the inputs $u_t$ to converge to $\eta_{t_i}$ at a similar rate as the states $x_t$ converge to $\theta_{t_i}$. When the algorithm receives a changed cost function, i.e., at time instances $t_i+1$, the second term on the right hand side of~\eqref{eq:bounded_inputs} is nonzero and allows $u_{t_i+1}$ to increase. Note that the bound depends on both $\eta$ and $\theta$, because the new cost function may be such that $\eta_{t_i} = \eta_{t_{i-1}}$, but $\theta_{t_i} \neq \theta_{t_{i-1}}$. In this case, $u_t$ has to be allowed to diverge from $\eta_{t_{i-1}}$ in order to reach the new optimal state $\theta_{t_i}$. 

Finally, we require the cost functions to be Lipschitz continuous, which is a common assumption in the literature (compare, e.g., \cite{Li2019,Nonhoff2022}).

\begin{assumption} \label{assump:Lipschitz} All cost functions $L_t\in\classCost$ are Lipschitz continuous, i.e., there exists $l>0$ such that
\[
\norm{L_t(u_1,x_1) - L_t(u_2,x_2)} \leq l \norm{(u_1,x_1) - (u_2,x_2)}
\]
for all $t \in \mathbb Z_{\geq0}$ and all $(u_1,x_1),(u_2,x_2) \in \mathbb R^m\times \mathbb R^n$.
\end{assumption}

Finally, we are ready to state a sufficient condition for asymptotic stability of the optimal steady state implying bounded regret for time-varying cost functions.

\begin{theorem} \label{thm:stab_regret}
	Let Assumptions~\ref{assump:posdef_cost} and~\ref{assump:avg_dwell_time}-\ref{assump:Lipschitz} hold and let $\varphi \geq \underline \varphi$ with $\underline \varphi$ from Lemma~\ref{lem:bounded_product}. The algorithm achieves bounded dynamic regret if there exists a constant $M$ such that
	\[
		\lim_{T\rightarrow \infty} \sum_{t=0}^T \sigma(t) \leq M. \numberthis \label{eq:upper_bound_sum_at}
	\]
\end{theorem}
\begin{proof}
	We prove the claim by deriving a regret upper bound assuming that~\eqref{eq:upper_bound_sum_at} holds. First, for any initial state $x_0\in\mathbb R^n$, Assumption~\ref{assump:Lipschitz} yields
	\begin{align*}
		\mathcal R_T(x_0,u) &= \sum_{t=0}^T L_t(u_t,x_t) - L_t(\eta_t,\theta_t) \\
		%&\leq l \sum_{t=0}^T \norm{(u_t,x_t) - (\eta_t,\theta_t)} \\
		&\leq l \sum_{t=0}^T \norm{x_t - \theta_t} + l\sum_{t=0}^T\norm{u_t-\eta_t}. \numberthis \label{eq:regret_bound_2sums}
	\end{align*}
	We proceed to bound the two sums in~\eqref{eq:regret_bound_2sums} separately. Since the cost function $L_t$ is revealed to the algorithm only at time instance $t+1$, i.e., with one time step delay, the algorithm stabilizes $\theta_{t_i}$ according to Assumption~\ref{assump:state_sequence} during the interval $t\in[t_i+1,t_{i+1}+1]$. Note that $x_{t_i+1}$ is both, the starting point of the sequence converging to $\theta_{t_i}$ and the endpoint of the sequence converging to $\theta_{t_{i-1}}$. Thus, we have for $i\in[0,N-1]$,
	\begin{align}
		&\sum_{t=t_i+1}^{t_{i+1}} \norm{x_t - \theta_t} \leq \norm{\theta_{t_{i+1}} - \theta_{t_{i}}} + \sum_{t=t_i+1}^{t_{i+1}} \norm{x_t - \theta_{t_{i}}} \label{eq:asymp_stab_proof_1}
	\end{align}
	and, by Assumption~\ref{assump:state_sequence},
	\begin{align*}
		&\sum_{t=t_i+1}^{t_{i+1}} \norm{x_t - \theta_{t_i}}	\leq k\norm{x_{t_i+1} - \theta_{t_i}} \sum_{t=t_i+1}^{t_{i+1}} \sigma(t-t_i-1) \\
		\refleq{\eqref{eq:upper_bound_sum_at}} &kM \Big(\norm{\theta_{t_i}-\theta_{t_{i-1}}} + \norm{x_{t_i+1} - \theta_{t_{i-1}}} \Big),
	\end{align*}
	Repeatedly applying Assumption~\ref{assump:state_sequence} yields
	\begin{align*}
		&\sum_{t=t_i+1}^{t_{i+1}} \norm{x_t - \theta_{t_i}} \\
		\leq &kM \norm{\theta_{t_i}-\theta_{t_{i-1}}} + kM\norm{x_{t_{i-1}+1} - \theta_{t_{i-1}}} k(\Delta \sigma)_i \\
		\begin{split}
			\leq &kM \norm{x_1-\theta_0} \prod_{j=1}^{i} k(\Delta \sigma)_j \\ 
			&\quad + kM \sum_{j=1}^{i} \left( \norm{\theta_{t_{j}}-\theta_{t_{j-1}}} \prod_{s=j}^{i-1} k(\Delta \sigma)_{s+1} \right).
		\end{split}
	\end{align*}
	Summing over $i\in[0,N-1]$ yields
	\begin{align*}
		&\sum_{i=0}^{N-1} \sum_{t=t_i+1}^{t_{i+1}} \norm{x_t - \theta_{t_i}} \leq kM \norm{x_1-\theta_0} \sum_{i=0}^{N-1} \left( \prod_{j=1}^{i} k(\Delta \sigma)_j \right)  \\
		&\hspace{10ex} + kM\sum_{i=0}^{N-1} \sum_{j=1}^{i} \norm{\theta_{t_{j}}-\theta_{t_{j-1}}} \left( \prod_{s=j}^{i-1} k(\Delta \sigma)_{s+1} \right).
	\end{align*}
	By using a case distinction on whether $N \geq N_0+1$ and defining $C_{N_0} = \sum_{i=0}^{N_0} k^i$, we get
	\begin{align*}
		&\sum_{i=0}^{N-1} \sum_{t=t_i+1}^{t_{i+1}} \norm{x_t - \theta_{t_i}}  \refleq{\eqref{eq:bound_sum_prod}} kM\max ( C_{N_0}, P) \norm{x_1-\theta_0} \\
		&\hspace{5ex} + kM\sum_{j=1}^{N-1} \left( \norm{\theta_{t_{j}}-\theta_{t_{j-1}}} \sum_{i=0}^{N-1-j} \prod_{s=1}^{i} k(\Delta \sigma)_{s+j} \right) \\
		\refleq{\eqref{eq:bound_sum_prod}}  &\hat C  \left( \norm{x_1-\theta_0} + \sum_{j=1}^{N-1} \norm{\theta_{t_{j}}-\theta_{t_{j-1}}} \right) \\
		\leq &\hat C \norm{x_1-\theta_0} + \hat C \sum_{t=1}^{T} \norm{\theta_{t}-\theta_{t-1}}, \numberthis \label{eq:asymp_stab_proof_2}
	\end{align*}
	where $\hat C = kM \max ( C_{N_0}, P)$. Combining the above results we get
	\begin{align*}
		&\sum_{t=0}^T \norm{x_t - \theta_t} = \norm{x_0-\theta_0} + \sum_{i=0}^{N-1} \sum_{t=t_i+1}^{t_{i+1}} \norm{x_t - \theta_t} \\
		\refleq{\eqref{eq:asymp_stab_proof_1}} &\norm{x_0-\theta_0} + \sum_{i=0}^{N-1} \norm{\theta_{t_{i+1}} - \theta_{t_{i}}} + \sum_{i=0}^{N-1} \sum_{t=t_i+1}^{t_{i+1}} \norm{x_t - \theta_{t_{i}}} \\
		\refleq{\eqref{eq:asymp_stab_proof_2}} &C_0 + \left( 1+ \hat C\right) \sum_{t=1}^{T} \norm{\theta_{t} - \theta_{t-1}}, \numberthis \label{eq:regret_bound_x-theta}
	\end{align*}
	where $C_0 = \norm{x_0-\theta_0} + \hat C \norm{x_1 - \theta_0}$. Note that $\norm{x_0-\theta_0}$ and $\norm{x_1 - \theta_0}$ are constants which only depend on the initialization $(u_0,x_0)$. It remains to bound $\sum_{t=0}^T \norm{u_t - \eta_t}$. Assumption~\ref{assump:bounded_inputs} yields
	\begin{align*}
		&\sum_{t=0}^T \norm{u_t-\eta_t} \leq \sum_{t=0}^T \norm{u_t-\eta_{t-1}} + \sum_{t=0}^T \norm{\eta_t - \eta_{t-1}} \\
		\refleq{\eqref{eq:bounded_inputs}} & k_u\sum_{t=0}^{T} \norm{u_{t-1} - \eta_{t-1}} + k_x \sum_{t=0}^T \norm{x_t - \theta_{t-1}} \\
		&+ k_\zeta \sum_{t=1}^T \norm{\zeta_{t-1} - \zeta_{t-2}} + \sum_{t=1}^T \norm{\eta_t - \eta_{t-1}}. 
	\end{align*}
	By Assumption~\ref{assump:bounded_inputs}, $u_{-1} = \eta_0$. Thus, rearranging yields
	\begin{align*}
		\begin{split}
		\sum_{t=0}^T \norm{u_t-\eta_t}\refleq{\eqref{eq:regret_bound_x-theta}} &\frac{k_x\left(2+\hat C\right) + k_\zeta}{1-k_u} \sum_{t=1}^T \norm{\theta_t - \theta_{t-1}} \\ +\frac{1+k_\zeta}{1-k_u} &\sum_{t=1}^T \norm{\eta_t-\eta_{t-1}} + \frac{k_x}{1-k_u} C_0
		\end{split} \numberthis \label{eq:regret_bound_u-eta}
	\end{align*}
	The desired regret bound then follows from inserting~\eqref{eq:regret_bound_x-theta} and~\eqref{eq:regret_bound_u-eta} into~\eqref{eq:regret_bound_2sums}. \hfill $\square$
\end{proof}

Theorem~\ref{thm:stab_regret} provides a sufficient condition for achieving bounded dynamic regret (for time-varying cost functions) if the closed loop is already known to be asymptotically stable (for any constant cost function). We note that the sufficient condition in Theorem~\ref{thm:stab_regret} and the necessary condition in Proposition~\ref{prop:necessary_condition} are similar, but not the same. More specifically, both conditions require the upper bound in Definition~\ref{def:asymp_stability} to be summable, but Theorem~\ref{thm:stab_regret} additionally assumes that the upper bound $\beta\in\mathcal{KL}$ is linear in $s = \norm{x_0-\theta}$. Thus, it is possible that there exists a $\beta\in\mathcal{KL^+}$ such that the optimal steady state is asymptotically stable according to Definition~\ref{def:asymp_stability}, thereby satisfying the necessary condition in Proposition~\ref{prop:necessary_condition}, but every linear bound $k\norm{x_0-\theta}\sigma(t)$, even if it exists, is not summable, i.e., does not satisfy the sufficient condition of Theorem~\ref{thm:stab_regret}. Combining Proposition~\ref{prop:necessary_condition} and Theorem~\ref{thm:stab_regret} in order to obtain a necessary and sufficient condition is an interesting direction for future research.

\begin{remark} (Algorithms that admit a state-space representation)
	If Assumption~\ref{assump:algorithm_state_space_rep} is satisfied, i.e., the Algorithm $\algorithm$ admits a state-space representation, some of the requirements for Theorem~\ref{thm:stab_regret} can be relaxed. More specifically, we assume that the extended state $(x,x^\algorithm)$ is asymptotically stable with respect to the set $\mathcal S^\theta_L = \{ (\theta,\theta^\algorithm) \}$ according to Assumption~\ref{assump:state_sequence}, where $\theta^\algorithm$ are the unique controller states $x^\algorithm$ such that $x^\algorithm\in\mathcal S_L$. Thereby, we allow the system states to diverge from $\theta$ initially if the algorithm states $x^\algorithm$ are not initialized correctly (compare Remark~\ref{rem:asymp_stab_general_algo}). Moreover, instead of Assumption~\ref{assump:bounded_inputs}, we require that $h_L^\algorithm(x^\algorithm,x)$ in Assumption~\ref{assump:algorithm_state_space_rep} is uniformly (in $L\in\classCost$) Lipschitz continuous, and that the variation of the optimal steady-state controller states $\theta_t^\algorithm$ is bounded by the variation of the optimal steady states and inputs, i.e., there exist constants $k_\theta^\algorithm,k_\eta^\algorithm>0$ such that for any two cost functions $L_1, L_2\in\classCost$, the corresponding steady states $(\eta_1,\theta_1),(\eta_2,\theta_2)$ and controller states $\theta^\algorithm_1,\theta^\algorithm_2$ satisfy $\norm{\theta_1^\algorithm-\theta_{2}^\algorithm} \leq k_\theta^\algorithm \norm{\theta_1-\theta_{2}} + k_\eta^\algorithm \norm{\eta_1-\eta_{2}}$. Then, we have bounded regret as in Theorem~\ref{thm:stab_regret} by similar arguments as above. Additionally, as discussed above Lemma~\ref{lem:bounded_product}, if there exists a compact set $\mathcal X \subset \mathbb R^n$ such that $x_t \in \mathcal X$ for all $t\in\mathbb Z_{\geq0}$ and if there exists $\kappa\in\mathcal{K_\infty}$, $\sigma\in\mathcal L$ such that $\norm{(x_t,x_t^\algorithm)}_{\mathcal S_L^\theta} \leq \kappa(\norm{(x_0,x_0^\algorithm)}_{\mathcal S_L^\theta})\sigma(t)$ for any constant cost function $L\in\classCost$ with $\kappa(s)$ Lipschitz at the origin, then Assumption~\ref{assump:state_sequence} is satisfied. Moreover, if the algorithm admits a state-space representation, then Lipschitz continuity at the origin of $\kappa(s)$ implies that the sufficient condition in Theorem~\ref{thm:stab_regret}, i.e.,~\eqref{eq:upper_bound_sum_at}, is satisfied as well \cite[Proposition~2]{Rakovic2019}.
\end{remark}

We close this section by two examples illustrating the application of Theorem~\ref{thm:stab_regret}. First, we continue Example~\ref{cex:stable->regret} and modify the proposed algorithm such that it achieves bounded regret. Thereafter, in Example~\ref{ex:exponential_stab}, we design an exponentially stabilizing algorithm and derive an average dwell time for this special case that ensures bounded dynamic regret.

\setcounter{example}{0}
\begin{example} (continued)
	Recall that we consider the \linebreak discrete-time integrator $x_{t+1} = x_t + u_t$. We design an improved version of our algorithm given by
	\[
		u_t = -\frac{2\tau+1}{(\tau+1)^2} x_t + \frac{2\tau+1}{(\tau+1)^2}\theta_{t-1},
	\]
	where $\theta_{-1} = x_0$ and $\tau$ as described above, such that, for a constant minimizer $\theta$, the closed loop achieves $\norm{x_t - \theta} = \frac{1}{t^2}\norm{x_0-\theta}$ for $t\in\mathbb Z_{\geq1}$ satisfying Assumption~\ref{assump:state_sequence} with \linebreak $\sigma(0) = 1$, $\sigma(t) = \frac{1}{t^2}$ for $t\in\mathbb Z_{\geq1}$, and $k=1$. Since $(\eta_t,\theta_t)$ are a steady state, we get $\eta_t=0$ for all $t$, which implies
	\[
		\norm{u_t - \eta_{t-1}} = \norm{u_t} = \frac{2\tau+1}{(\tau+1)^2}\norm{x_t - \theta_{t-1}},
	\]
	i.e., Assumption~\ref{assump:bounded_inputs} is satisfied with $k_\zeta=k_u=0$ and $k_x = 1$ since $\tau \geq 0$. Thus, the improved algorithm achieves bounded regret if the cost functions satisfy the remaining assumptions of Theorem~\ref{thm:stab_regret}, i.e., Assumption~\ref{assump:avg_dwell_time} with $\varphi\geq \underline\varphi$ and Assumption~\ref{assump:Lipschitz}. In order to ensure the former condition, we show that we can choose any $\underline \varphi>1$, i.e.,~\eqref{eq:bound_sum_prod} holds for any $\underline \varphi > 1$ and $N_0 > 0$. To this end, fix $N''$ to be the smallest integer greater than or equal to $(\underline\varphi N_0+1)/(\underline\varphi-1)$ which implies $t_{N''} > \underline \varphi (N''-N_0) \geq N''+1$, where the first inequality follows from Assumption~\ref{assump:avg_dwell_time}. Then, we get
	\[
		\prod_{j=1}^{N''} k(\Delta \sigma)_j \leq \frac{1}{4}
	\]
	since $(\Delta \sigma)_j \leq \frac{1}{4}$ for at least one $j \in [1,N'']$. For any $N'\in\mathbb Z_{>N_0}$ we have $N' \leq rN''$ for some $r\in\mathbb Z_{>0}$ and
	\[
		\sum_{i=0}^{N'} \prod_{j=1}^{i} (\Delta \sigma)_j \leq \sum_{i=0}^{rN''} \prod_{j=1}^{i} (\Delta \sigma)_j \leq N'' \sum_{i=0}^{r} \left(\frac{1}{4}\right)^i \leq \frac{4N''}{3},
	\]
	i.e., \eqref{eq:bound_sum_prod} holds. Thus, if the cost functions do not change (on average) at every time step and satisfy the remaining assumptions in Theorem~\ref{thm:stab_regret}, the improved algorithm achieves bounded dynamic regret.\hfill $\square$
\end{example}

The last example illustrates how exponential stability implies bounded regret.

\begin{example} \label{ex:exponential_stab}
	Consider a linear system $x_{t+1} = Ax_t + Bu_t$. We design an algorithm~$\algorithm$ that, at each time instance~$t$, solves an optimization problem to obtain $(\eta_{t-1},\theta_{t-1})$ and applies
	\[
		u_t = K(x_t - \theta_{t-1}) + \eta_{t-1},
	\]
	where $K$ is chosen such that the closed-loop dynamics \linebreak $(A+BK)$ are Schur stable. Then, for a constant cost function~$L\in\classCost$, we have $\theta = A\theta + B\eta$ by Assumption~\ref{assump:posdef_cost} and
	\begin{align*}
		\norm{x_t - \theta} &= \norm{A(x_{t-1}-\theta) + B(u_{t-1}-\eta)} \\
		%&= \norm{(A+BK)(x_{t-1} - \theta)} \\
		&\leq \norm{(A+BK)^t}\norm{x_0-\theta}.
	\end{align*}
	Since $(A+BK)$ is Schur stable, there exist constants $c\geq1$ and $\lambda\in(0,1)$ such that
	\[
		\norm{x_t - \theta} \leq c\lambda^t \norm{x_0-\theta},
	\]
	i.e., the closed loop is exponentially stable with respect to~$\theta$ and Assumption~\ref{assump:state_sequence} is satisfied with $\sigma(t) = \lambda^t$ and $k=c$. Moreover, since
	\[
		\norm{u_t - \eta_{t-1}} = \norm{K(x_t - \theta_{t-1})} \leq \norm{K} \norm{x_t - \theta_{t-1}},
	\]
	Assumption~\ref{assump:bounded_inputs} is satisfied with $k_x = \norm{K}$ and $k_\zeta=k_u = 0$. Finally, we derive an explicit expression for $\underline \varphi$ to ensure that~\eqref{eq:bound_sum_prod} is satisfied. To this end, note that Lemma~\ref{lem:bounded_product} holds for $\underline \varphi = -\frac{\ln(k)}{\ln(\lambda)}+\varphi_0>0$ and any $\varphi_0>0$, since we get $t_i > \underline\varphi(i-N_0)$ by Assumption~\ref{assump:avg_dwell_time} and, therefore, 
	\begin{align*}
		&\sum_{i=0}^{N'} \prod_{j=1}^{i} k(\Delta \sigma)_j = \sum_{i=0}^{N'} k^{i} \lambda^{t_i} \\
		< &\sum_{i=0}^{N'} e^{\ln(k) i} e^{-\ln(k)(i-N_0)} \lambda^{\varphi_0(i-N_0)}  \\ 
		= &k^{N_0} \lambda^{-\varphi_0N_0} \sum_{i=0}^{N'} \lambda^{\varphi_0 i}
		\leq \frac{ k^{N_0} \lambda^{-\varphi_0N_0} }{1-\lambda^{\varphi_0}},
	\end{align*}
	i.e., \eqref{eq:bound_sum_prod} holds. Thus, if the cost functions satisfy the remaining assumptions of Theorem~\ref{thm:stab_regret}, then the algorithm~$\mathcal A$ achieves bounded dynamic regret.
\end{example}

%%%%%%%%%%%%%%%%%%%%%%% Conclusion %%%%%%%%%%%%%%%%%%%%%%%%%%%%%%5
	
\section{CONCLUSION} \label{sec:conclusion}

In this paper, we study the relation between dynamic regret and asymptotic stability. Loosely speaking, our results suggest that achieving bounded dynamic regret is a stronger property than asymptotic stability, because \linebreak bounded dynamic regret requires a certain minimal rate of convergence, whereas asymptotic stability only characterizes the asymptotic behavior of the closed loop. In particular, we show that, under rather mild assumptions, bounded dynamic regret implies asymptotic stability (compare Theorem~\ref{thm:regret_asymp_stab}), whereas asymptotic stability does, in general, not imply bounded dynamic regret (compare Example~\ref{cex:stable->regret} and Proposition~\ref{prop:necessary_condition}). In order to derive sufficient conditions for the case of bounded regret implying asymptotic stability, however, we require additional, more restrictive assumptions (compare Theorem~\ref{thm:stab_regret}). Therefore, future works includes relaxing these additional assumptions. In particular, we conjecture that Assumption~\ref{assump:avg_dwell_time} could be relaxed as discussed above. Moreover, one could consider local asymptotic stability and constraints instead of the global results provided in this paper. Another interesting possibility for future work is considering other common applications of dynamic regret, e.g., with respect to a priori unknown process noise or unknown system dynamics instead of time-varying cost functions. Finally, future work includes studying the relation between asymptotic stability and dynamic regret for benchmarks that do not achieve bounded dynamic regret with respect to the optimal steady state themselves (compare Remark~\ref{rem:different_benchmark}).

%%%%%%%%%%%%% Appendix %%%%%%%%%%%%%%%%%%%
	
\appendix

\section{Proof of Lemma~\ref{lem:summable_KL_function}}

If there exists $\beta\in\mathcal{KL^+}$ such that~\eqref{eq:proposition} holds, then we have for any constant cost function $L\in\classCost$, initial state $x_0\in\mathbb R^n$, and $s=\norm{x_0-\theta}$
\[
	\lim_{T\rightarrow\infty} \sum_{t=0}^T \norm{x_t - \theta} \refleq{\eqref{eq:proposition}} \lim_{T\rightarrow\infty} \sum_{t=0}^T \beta(s,t) \leq B(s) =: D(s),
\]
where the second inequality follows from $\beta\in\mathcal{KL^+}$. 

For the converse implication, fix any constant cost function~$L{\in\classCost}$ and assume that for any initial state $x_0\in\mathbb R^n$ there exists a finite constant $D(\norm{x_0-\theta})\geq 0$ such that
\begin{equation} \label{eq:summable_traj}
	\lim_{T\rightarrow\infty} \sum_{t=0}^T \norm{x_t - \theta} \leq D(\norm{x_0-\theta}) < \infty.
\end{equation}
Moreover, there exists $\bar\beta\in\mathcal{KL}$ (but not necessarily in $\mathcal{KL^+}$) such that~\eqref{eq:proposition} holds because the closed loop is asymptotically stable. Then, we construct $\beta\in\mathcal{KL^+}$ such that~\eqref{eq:proposition} is satisfied. 

To this end, for any $s\geq0$ and $t\in\mathbb Z_{\geq0}$, consider $\phi(s,t)$ as defined in~\eqref{eq:def_phi}. Note that $\phi(0,t) = 0$, $\phi(\cdot,t)$ is nondecreasing for any fixed~$t$, and $\phi(\cdot,t)$ is upper bounded by $\bar\beta(\cdot,t) \in \mathcal K$ because the closed loop is asymptotically stable by the assumptions of Lemma~\ref{lem:summable_KL_function}. Hence, for every $t\in\mathbb Z_{\geq0}$, $\phi(\cdot,t)$ is continuous at the origin. Therefore, for any $\epsilon>0$, there exists a function $\phi_c:\mathbb R_{[0,\epsilon]}\times\mathbb Z_{\geq0} \rightarrow \mathbb R_{\geq0}$ such that, for any fixed $t\in\mathbb{Z}_{\geq0}$, $\phi_c(\cdot,t)$ is continuous, nondecreasing, $\phi_c(0,t) = 0$, $\phi_c(\epsilon,t) = \phi(2\epsilon,t)$, $\phi_c(s,t) \geq \phi(s,t)$ for all $s\in[0,\epsilon]$, and $\phi_c(s,\cdot)$ is nonincreasing for any fixed $s\in[0,\epsilon]$. Next, we define the function
\[
	\Phi(s,t) {=} \! \begin{cases} \phi_c(s,t) &\!\text{if } s\in[0,\epsilon] \\ \phi(i\epsilon,t) + \frac{s-(i{-}1)\epsilon}{\epsilon} \Delta\phi(i\epsilon,t) &\!\text{if } s \in ((i{-}1)\epsilon,i\epsilon]\\&\!\phantom{\text{if }}i = 2,3,\dots \end{cases}\!,
\]
where $\Delta\phi(i\epsilon,t) = \phi((i+1)\epsilon,t)-\phi(i\epsilon,t)$. The function~$\Phi$ is piece-wise linear for any fixed $t\in\mathbb Z_{\geq0}$ and $s\geq\epsilon$ and it holds that $\phi(s+2\epsilon,t)\geq\Phi(s,t) \geq \phi(s,t)$ for any $s\geq0$ and $t\in\mathbb Z_{\geq0}$. In the following, we will show that for any fixed $s\geq0$, $\Phi(s,\cdot)\in\mathcal{L}$, and for any fixed $t\in\mathbb{Z}_{\geq0}$, $\Phi(\cdot,t)$ is nondecreasing, continuous and $\Phi(0,t) = 0$. First, fix any $s\geq0$. Then, $\Phi(s,\cdot)$ is nonincreasing because both $\phi(s,\cdot)$ and $\phi_c(s,\cdot)$ are nonincreasing by definition. Moreover,
\[
	0 \leq \lim_{t\rightarrow\infty}\Phi(s,t) \leq \lim_{t\rightarrow\infty} \phi(s+2\epsilon,t) \leq \lim_{t\rightarrow\infty} \bar\beta(s+2\epsilon,t) = 0,
\]
i.e., $\Phi(s,\cdot)\in\mathcal L$. Second, $\Phi(\cdot,t)$ is continuous by definition for any fixed $t\in\mathbb Z_{\geq0}$ and $\Phi(0,t) = 0$ because $\phi_c(0,t)=0$. Additionally, $\Phi(\cdot,t)$ is nondecreasing since both $\phi_c(\cdot,t)$ and $\phi(\cdot,t)$ are nondecreasing.

Finally, choose any function $\beta'\in\mathcal{KL^+}$ which satisfies, for any $s\geq0$, $2\beta'(s,t) \geq \bar \beta(s,t)$ for all $t< \underline T$ and $\beta'(s,\underline T-1) \geq \Phi(s,\underline T)$, where $\underline T$ is from Assumption~\ref{assump:bounded_convergence}. Then, we define
\[
	\beta(s,t) = \begin{cases} 
		2\beta'(s,t)&\text{if } t < \underline T \\
		\beta'(s,t) + \Phi(s,t) &\text{if } t\geq \underline T
		\end{cases}.
\]	
We claim that $\beta$ is the desired function, i.e., $\beta\in\mathcal{KL^+}$ and it satisfies~\eqref{eq:proposition}. Note that $\beta\in\mathcal{KL}$ because $\beta'\in\mathcal{KL}$, $\Phi(s,\underline T) \leq \beta'(s,\underline T-1)$ for any $s$, and by the properties of $\Phi$ discussed above. We first show that \eqref{eq:proposition} is satisfied for any initial state $x_0 \in \mathbb R^n$. Let $s=\norm{x_{0}-\theta}$. We have
\begin{align*}
	\norm{x_{t} - \theta} &\leq \bar \beta(s,t) \leq 2\beta'(s,t) = \beta(s,t)
	\intertext{for $t < \underline T$ and for $t\geq\underline T$ we get}
	\norm{x_{t}-\theta} &\leq \phi(s,t) \leq \Phi(s,t) \leq \beta(s,t).
\end{align*}
Second, we have for any $s\geq0$
\begin{align*}
	\lim_{T\rightarrow\infty} \sum_{t=0}^T \beta(s,t) &= C_{\underline T} + \lim_{T\rightarrow\infty} \sum_{t=\underline T}^T \Phi (s,t) \\
	&\leq C_{\underline T} + \lim_{T\rightarrow\infty} \sum_{t=\underline T}^T \phi(s+2\epsilon,t),
\end{align*}
where $C_{\underline T} = \sum_{t=0}^{\underline T-1} \beta'(s,t) + \lim_{T\rightarrow \infty} \sum_{t=0}^T \beta'(s,t)$ is finite because $\beta'\in\mathcal{KL^+}$. Finally, by Assumption~\ref{assump:bounded_convergence} there exists a closed-loop trajectory $\bar x = \{\bar x_t\}_{t=0}^\infty$ such that we get
\begin{align*}
	\lim_{T\rightarrow\infty} \sum_{t=0}^T \beta(s,t) &\leq C_{\underline T} + \mu \lim_{T\rightarrow\infty} \sum_{t=\underline T}^T \norm{\bar x_t - \theta} \\
	&\refleq{\eqref{eq:summable_traj}} C_{\underline T} + \mu D(\norm{\bar x_0 - \theta}) =: B(s),
\end{align*}
i.e., $\beta\in\mathcal{KL^+}$.\hfill $\square$

\section{Proof of Lemma~\ref{lem:bounded_product}}

First, we claim that for any $N_0 \in \mathbb Z_{\geq0}$ there exists $\underline \varphi$ and $0<\delta<1$ such that
\[
	\prod_{j=1}^{N''} k(\Delta \sigma)_j \leq \delta \numberthis \label{eq:prod_bound_delta}
\]
for any $N''\in[N_0+1,2N_0+1]$ and $\varphi \geq \underline\varphi$. Let 
\[
	\bar \delta(N'') := \min_{\{t_i\}_{i=0}^{N''}} \prod_{j=1}^{N''} k(\Delta \sigma)_j \text{ s.t.~\eqref{eq:avg_dwell_time}} 
\]
be the smallest upper bound of the above product for any switching sequence $\{t_i\}_{i=0}^{N''}$ that satisfies Assumption~\ref{assump:avg_dwell_time}, which exists and satisfies $\bar \delta(N'')>0$ due to \mbox{$\sigma\in\mathcal L$}. Since $\sigma(\cdot)$ is nonincreasing, $\bar \delta(N'')$ is nonincreasing with respect to $\varphi$. Moreover, since $\lim_{t\rightarrow\infty} \sigma(t) = 0$ by Assumption~\ref{assump:state_sequence}, there exists $\underline \varphi$ such that $\bar \delta < 1$ holds for any $N''\in[N_0+1,2N_0+1]$, i.e., \eqref{eq:prod_bound_delta}. Next, we bound the sum~\eqref{eq:bound_sum_prod} for any $N'\in\mathbb I_{>N_0}$. Note that \mbox{$N' = r(N_0+1)+N''$} for some $r\in\mathbb Z_{\geq0}$ and $N''\in[N_0+1,2N_0+1]$. Thus, \mbox{$1 \geq (\Delta \sigma)_t > 0$} and $k>1$ lead to
\begin{align*}
	&\sum_{i=0}^{N'} \prod_{j=1}^{i} k(\Delta \sigma)_j = \sum_{i=0}^{r(N_0+1)+N''} \prod_{j=1}^{i} k(\Delta \sigma)_j \\
	= &\sum_{i=0}^{N''-1} \prod_{j=1}^{i} k(\Delta \sigma)_j + \sum_{i=N''}^{r(N_0+1)+N''} \prod_{j=1}^{i} k(\Delta \sigma)_j \\
	\leq &\sum_{i=0}^{2N_0} k^i + \sum_{i=N_0+1}^{(r+2)(N_0+1)} \prod_{j=1}^{i} k(\Delta \sigma)_j \numberthis \label{eq:Summenformel} \\
	\refleq{\eqref{eq:prod_bound_delta}} &\frac{1-k^{2N_0+1}}{1-k} + \sum_{i=1}^{r+2} (N_0+1)\delta^i \\
	\leq &\frac{1-k^{2N_0+1}}{1-k} + \frac{N_0+1}{1-\delta} =: P.
\end{align*}
If $k=1$, the first sum in~\eqref{eq:Summenformel} is bounded by $2N_0+1$. \hfill $\square$

\bibliographystyle{elsarticle-num}
\bibliography{refs}

\end{document}